\input amstex
\documentstyle{amsppt}
\magnification=\magstep1
 \hsize 13cm \vsize 18.35cm \pageno=1
\loadbold \loadmsam
    \loadmsbm
    \UseAMSsymbols
\topmatter
\NoRunningHeads
\title A study on the $q$-Euler numbers and the fermionic $q$-integral of the product of several type $q$-Bernstein polynomials on $\Bbb Z_p$
 \endtitle
\author
  Taekyun Kim \endauthor

\abstract
In this paper, we investigate some properties between $q$-Bernstein polynomials and $q$-Euler numbers by using fermionic $p$-adic $q$-integrals on $\Bbb Z_p$. From these properties, we derive some interesting identities related to  the $q$-Euler numbers.
\endabstract
\thanks  2000 AMS Subject Classification: 11B68, 11S80, 60C05, 05A30
\newline  The present Research has been conducted by the research
Grant of Kwangwoon University in 2010
\endthanks
\endtopmatter

\document

{\bf\centerline {\S 1. Introduction}}

 \vskip 20pt

Let $p$ be a fixed odd prime number. Throughout this paper $\Bbb Z_p$, $\Bbb Q_p$, and $\Bbb C_p$ will denote the ring of $p$-adic rational integers, the fields of $p$-adic rational numbers, and the completion of algebraic closure of $\Bbb Q_p$, respectively. Let $v_p$ be the normalized exponential
valuation of $\Bbb C_p$ with $|p|_p=p^{-v_p(p)}=\frac{1}{p}$. When one talks of $q$-extension, $q$ is variously considered as an indeterminate, a complex number $q\in \Bbb C$ or $p$-adic number $q\in \Bbb C_p$. If $q\in \Bbb C$, one normally assumes $|q|<1$, and if $q\in \Bbb C_p$, one normally assumes $|1-q|_p<1.$ We use the notations of $q$-numbers as $[x]_q=\frac{1-q^x}{1-q}$ and $[x]_{-q}=\frac{1-(-q)^x}{1+q}.$
Let $C[0, 1]$ denote the set of continuous function on $[0, 1]$. For $f\in C[0, 1]$, $q$-Bernstein operator is introduced in [2]
as follows:
$$ \Bbb B_{n, q}(f|x)=\sum_{k=0}^n f(\frac{k}{n})\binom{n}{k}[x]_q^k[1-x]_{\frac{1}{q}}^{n-k}=
\sum_{k=0}^n f(\frac{k}{n})B_{k,n}(x, q). \tag1 $$
Here $\Bbb B_{n,q}(f|x)$ is called $q$-Bernstein operator of order $n$ for $f$. For $k, n \in \Bbb Z_{+}(=\Bbb N \cup \{0\})$, the $q$-Bernstein polynomials of degree $n$ is defined by
$$B_{k,n}(x,q)=\binom{n}{k}[x]_q^k[1-x]_{\frac{1}{q}}^{n-k}, \text{ (see [2]).}\tag2$$
Let $UD(\Bbb Z_p)$ be the space of uniformly differentiable functions on $\Bbb Z_p$. For $f\in UD(\Bbb Z_p)$, the fermionic $p$-adic $q$-integrals on $\Bbb Z_p$ is defined by
$$I_q(f)=\int_{\Bbb Z_p}f(x) d\mu_{-q}(x)=\lim_{N\rightarrow \infty}\frac{1}{[p^N]_{-q}} \sum_{x=0}^{p^N-1}f(x)(-q)^x, \text{ (see [3, 4, 5])}.\tag3$$
 As well known definition, Euler polynomials are defined by
$$ \frac{2}{e^t+1}e^{xt}=\sum_{n=0}^{\infty}E_n(x)\frac{t^n}{n!}, \text{ (see [1-11]).}\tag4$$
In the special case, $x=0$, $E_n(0)=E_n$ are called the $n$-th Euler numbers. By (4), we see that the recurrence formula of Euler numbers is given by
$$E_0=1, \text{ and }  (E+1)^n+E_n=0 ,  \text{  if  $n>0$, (see [2, 4, 5]),}\tag5$$
with the usual convention of replacing  $E^n$ by $E_n$.
When one talks of $q$-analogue, $q$ is variously considered as an indterminate, a complex number $q\in \Bbb C$, or a $p$-adic
number $q\in \Bbb C_p$. If $q\in \Bbb C$, we normally assume $|q|<1$. If $q\in \Bbb C_p$, we normally always assume that $|1-q|_p<1$.
As the $q$-extension of (5), the $q$-Euler numbers are defined  by
$$ \xi_{0,q}=1, \text{ and } q(q\xi+1)^n+\xi_{n,q}=0, \text{ if $ n>0,$ (see [5]),} \tag6$$
 with the usual convention of replacing $\xi^n$ by $\xi_{n,q}$. The $q$-Euler polynomials are also defined by
 $$\xi_{n,q}(x)=\sum_{l=0}^n\binom{n}{l}[x]_q^{n-l}q^{lx}\xi_{l, q},  \text{ for $n\in\Bbb Z_{+},$ (see [5])}.\tag7$$
 By (6) and (7), we easily get
 $$q\xi_{n,q}(1)+\xi_{n,q}=1,  \text{ if $n>0$, (see [4, 5]).}\tag8$$
 In [5], Witt's formula for the $q$-Eula polynomial is given by
 $$\xi_{n,q}(x)=\int_{\Bbb Z_p}[y+x]_q^n d\mu_{-1}(y)=\frac{[2]_q}{(1-q)^n}\sum_{l=0}^n\binom{n}{l}(-1)^l \frac{q^{lx}}{
1+q^l}, \text{ for $ n\geq 0.$}\tag9$$
By (9), we eaily see that
$$\xi_{n,q}(x)=\sum_{l=0}^n\binom{n}{l}q^{lx}\xi_{l,q}=(q^x\xi+1)^n.$$
In this paper, we investigate some properties between $q$-Bernstein polynomials and $q$-Euler numbers by using fermionic $p$-adic $q$-integrals on $\Bbb Z_p$. From these properties, we derive some interesting identities related to  the $q$-Euler numbers.

\vskip 20pt

{\bf\centerline {\S 2. $q$-Euler numbers and  $q$-Bernstein  Polynomials}} \vskip 10pt

  In this section we assume that $q\in \Bbb C_p$ with $|1-q|_p<1.$ From (6), (7) and (9), we note that
  $$\aligned
  q^2\xi_{k,q}(2)-q-q^2&=q^2\left(q(q\xi +1)+1 \right)^k-q-q^2\\
  &=q^2\left(1+\sum_{l=1}^k\binom{k}{l}q^l(q\xi+1)^l\right)-q-q^2\\
  &=-q\sum_{l=1}^k\binom{k}{l}q^l\xi_{l,q}-q=-q(q\xi+1)^k=\xi_{k, q}, \text{ if $k>0$.}
  \endaligned\tag10$$
Therefore, we obtain the following  theorem.

\proclaim{ Theorem 1}
For $n\in \Bbb N$, we have
$$\xi_{n,q}(2)=1+\frac{1}{q}+\frac{1}{q^2}\xi_{n,q}.$$
\endproclaim

  By (9), we easily get
  $$\int_{\Bbb Z_p}[x]_q^n d\mu_{-q}(x)=\xi_{n,q}(0)=\xi_{n,q}=\frac{[2]_q}{(1-q)^n}\sum_{l=0}^n\binom{n}{l}(-1)^l \frac{1}{1+q^l}, \text{ for $n\geq 0$}, $$
  and
  $$\xi_{n,q}(x)=\int_{\Bbb Z_p}[x+y]_q^n d\mu_{-q}(y)=\sum_{l=0}^n \binom{n}{l}q^{lx}\int_{\Bbb Z_p}[y]_q^l d\mu_{-q}(y)=
  \sum_{l=0}^n \binom{n}{l}[x]_q^{n-l}q^{lx}\xi_{l,q}.\tag11$$
It is easy to show that
$[1-x]_{\frac{1}{q}}^n=(-1)^nq^n[x-1]_q^n.$ Thus, we have
$$\int_{\Bbb Z_p} [1-x]_{\frac{1}{q}}^n d\mu_{-q}(x)=(-1)^nq^n\int_{\Bbb Z_{p}}[x-1]_q^nd\mu_{-q}(x)=(-1)^nq^n \xi_{n,q}(-1). \tag12$$

  By (11) and (12), we obtain the following proposition.
 \proclaim{ Proposition 2}
 For $n\in\Bbb Z_{+}$, we have
 $$\int_{\Bbb Z_p}[1-x]_{\frac{1}{q}}^nd\mu_{-q}(x)=\sum_{l=0}^n\binom{n}{l}(-1)^l\xi_{l,q}=(-1)^nq^n \xi_{n,q}(-1). $$
 \endproclaim
From (3), we can derive the following equation (13).
$$\int_{\Bbb Z_p}[1-x+x_1]_{\frac{1}{q}}^n d\mu_{-\frac{1}{q}}(x_1)=(-1)^n q^n \frac{[2]_q}{(1-q)^n}\sum_{l=0}^n\binom{n}{l}
(-1)^l \frac{q^{lx}}{1+q^{l+1}}.\tag13$$
By (9) and (13), we obtain the following proposition.

\proclaim{ Proposition 3}
For $n\in \Bbb Z_{+}$, we have
$$\int_{\Bbb Z_p}[1-x+x_1]_{\frac{1}{q}}^n d\mu_{-\frac{1}{q}}(x_1)=(-1)^n q^n \int_{\Bbb Z_p}[x+x_1]_q d\mu_{-q}(x_1). \tag14$$
\endproclaim

By (12) and (14), we get
$$q^n(-1)^n\xi_{n,q}(-1)=(-1)^nq^n\int_{\Bbb Z_p}[x_1-1]_q^n d\mu_{-q}(x_1) =\int_{\Bbb Z_p}[2+x_1]_{\frac{1}{q}}^nd\mu_{-\frac{1}{q}}(x_1)
=\xi_{n, \frac{1}{q}}(2).\tag15$$

Therefore, by (15),  we obtain the following corollary.
\proclaim{ Corollary 4}
For $n\geq 0$, we have
$$\xi_{n, \frac{1}{q}}(2)=(-1)^nq^n \xi_{n,q}(-1).$$
\endproclaim
 From (12) and Corollary 4, we have
 $$\int_{\Bbb Z_p}[1-x]_{\frac{1}{q}}^n d\mu_{-q}(x)=(-1)^nq^n\xi_{n,q}(-1)=\xi_{n, \frac{1}{q}}(2).\tag16$$
By (12), (16) and Theorem 1, we obtain the following theorem.
\proclaim{Theorem 5}
For $n \in \Bbb N$, we have
$$\int_{\Bbb Z_p}[1-x]_{\frac{1}{q}}^n d\mu_{-q}(x)=q^2 \int_{\Bbb Z_p}[x]_{\frac{1}{q}}^n d\mu_{-\frac{1}{q}}(x)+q+1
=[2]_q+q^2\xi_{n, \frac{1}{q}}.$$
\endproclaim

Now we consider  the $p$-adic analogue of (1) and (2) on $\Bbb Z_p$. For $x\in \Bbb Z_p$ and $f\in UD(\Bbb Z_p)$, we define
$p$-adic $q$-Bernstein operator of order $n$ for $f$ as follows:
$$\Bbb B_{n, q}(f|x)=\sum_{k=0}^n f(\frac{k}{n})\binom{n}{k}[x]_q^k[1-x]_{\frac{1}{q}}^{n-k}=\sum_{k=0}^nf(\frac{k}{n})B_{k,n}(x,q).$$
For $n, k \in \Bbb Z_{+}$, the $p$-adic $q$-Bernstein polynomials of degree $n$ is defined by
$$B_{k, n}(x, q)=\binom{n}{k}[x]_q^k[1-x]_{\frac{1}{q}}^{n-k}, \text{ where $ x\in \Bbb Z_{p}$}. \tag17$$
Taking the fermionic $p$-adic $q$-integral on $\Bbb Z_p$ for one $q$-Bernstein polynomials in (17), we get
$$\aligned
\int_{\Bbb Z_p}B_{k, n}(x, q)d\mu_{-q}(x)&=\binom{n}{k}\int_{\Bbb Z_p}[x]_q^k[1-x]_{\frac{1}{q}}^{n-k} d\mu_{-q}(x)\\
&=\binom{n}{k}\sum_{l=0}^{n-k}\binom{n-k}{l}(-1)^l\int_{\Bbb Z_{p}}[x]_q^{k+l}d\mu_{-q}(x)\\
&=\binom{n}{k}\sum_{l=0}^{n-k}\binom{n-k}{l}(-1)^l \xi_{k+l, q}.
\endaligned\tag18$$
From (17), we note that
$$B_{k, n}(x, q)=B_{n-k, n}(1-x, \frac{1}{q}), \text{ where $n, k \in \Bbb Z_{+}$ and $ x\in \Bbb Z_p$}. \tag19$$
By (19), we get
$$\aligned
\int_{\Bbb Z_p}B_{k, n}(x, q)d\mu_{-q}(x)&=\int_{\Bbb Z_p}B_{n-k, n}(1-x, \frac{1}{q})d\mu_{-q}(x)\\
&=\binom{n}{k}\sum_{l=0}^{k}\binom{k}{l}(-1)^{k+l}\int_{\Bbb Z_p}[1-x]_{\frac{1}{q}}^{n-l}d\mu_{-q}(x)\\
&=\binom{n}{k}\sum_{l=0}^k \binom{k}{l}(-1)^{k+l}\xi_{n-l, \frac{1}{q}}(2).
\endaligned\tag20$$
By comparing the coefficients on the both sides of (18) and (20), we obtain the following theorem.

\proclaim{Theorem 6}
For $n, k \in \Bbb Z_{+}$ with $n\geq k$, we have
$$\sum_{l=0}^{n-k}\binom{n-k}{l}(-1)^l \xi_{k+l,q}=\sum_{l=0}^k \binom{k}{l}(-1)^{k+l}\xi_{n-l, \frac{1}{q}}(2). $$
\endproclaim
By Theorem 1 and Theorem 6, we obtain the following corollary.

\proclaim{Corollary 7}
Let $n, k \in \Bbb Z_{+}$ with $n>k.$ Then we have
$$\sum_{l=0}^{n-k}\binom{n-k}{l}(-1)^l \xi_{k+l, q}=q^2\sum_{l=0}^k \binom{k}{l}(-1)^{k+l}\xi_{n-l, \frac{1}{q}}, \text{ if $k >0$.}$$
Moreover,
$$\sum_{l=0}^n \binom{n}{l}(-1)^l \xi_{l, q}=[2]_q+ q^2\xi_{n, \frac{1}{q}}. $$
\endproclaim
For $ m, n, k \in \Bbb Z_+ $ with $m, n > k$, we consider the fermionic $p$-adic $q$-integral for multiplication of two $q$-Bernstein polynomials
on $\Bbb Z_p$.
$$\aligned
&\int_{\Bbb Z_p}B_{k, n}(x,q)B_{k, m}(x, q)d\mu_{-q}(x)=\binom{n}{k}\binom{m}{k}\int_{\Bbb Z_{p}}[x]_q^{2k}[1-x]_{\frac{1}{q}}^{n+m-2k}d\mu_{-q}(x)\\
&=\binom{n}{k}\binom{m}{k}\sum_{l=0}^{2k}\binom{2k}{l}(-1)^{l+2k}\int_{\Bbb Z_p}[1-x]_{\frac{1}{q}}^{n+m-l}d\mu_{-q}(x)\\
&=\binom{n}{k}\binom{m}{k}\sum_{l=0}^{2k}\binom{2k}{l}(-1)^{l+2k}\xi_{n+m-l, \frac{1}{q}}(2).
\endaligned\tag21$$
From (21) and Theorem 1, we can derive the following equation (22).
$$\aligned
&\int_{\Bbb Z_p}B_{k,n}(x,q)B_{k, m}(x,q)d\mu_{-q}(x)=\binom{n}{k}\binom{m}{k}\sum_{l=0}^{2k}\binom{2k}{l}(-1)^{l+2k}\xi_{n+m-l, \frac{1}{q}}(2)\\
&=\binom{n}{k}\binom{m}{k}\sum_{l=0}^{2k}\binom{2k}{l}(-1)^{l+2k}\left([2]_q+q^2\xi_{n+m-l, \frac{1}{q}}\right).
\endaligned\tag22$$
For $n, m, k \in \Bbb Z_{+}$, it is easy to show that
$$\int_{\Bbb Z_p}B_{k, n}(x, q)B_{k, m}(x, q) d\mu_{-q}(x)
=\binom{n}{k}\binom{m}{k}\sum_{l=0}^{n+m-2k}\binom{n+m-2k}{l}(-1)^l \xi_{l+2k, q}. \tag23$$
Therefore, by (22) and (23), we obtain the following theorem.
\proclaim{ Theorem 8}
Let $m, n, k \in \Bbb Z_{+}$ with $m+n>2k$. Then we have
$$\sum_{l=0}^{m+n-2k}\binom{n+m-2k}{l}(-1)^l\xi_{l+2k,q}=q^2\sum_{l=0}^{2k}\binom{2k}{l}(-1)^{l+2k}\xi_{n+m-l, \frac{1}{q}}, \text{ if $k >0$}.$$
Moreover,
$$\sum_{l=0}^{m+n}\binom{n+m}{l}(-1)^l \xi_{l, q} =q^2\xi_{n+m, \frac{1}{q}}+ [2]_q .$$
\endproclaim
By induction hypothesis, we obtain the following theorem.
 \proclaim{Theorem 9}
 Let $n_1, n_2, \cdots, n_s, k \in \Bbb Z_{+} $ $(s\in \Bbb N)$ with $n_1+n_2+\cdots+n_s>sk$. Then we have
 $$\aligned
 &\sum_{l=0}^{n_1+\cdots+n_s-sk}\binom{n_1+\cdots+n_s-sk}{l}(-1)^l\xi_{l+sk,q}\\
 &= q^2\sum_{l=0}^{sk}\binom{sk}{l}(-1)^{sk-l}\xi_{n_1+\cdots+n_s-l, \frac{1}{q}},
 \text{ if $k>0$.}\endaligned$$
 Moreover,
 $$\sum_{l=0}^{n_1+\cdots+n_s} \binom{n_1+\cdots+n_s}{l}(-1)^l\xi_{l, q}=[2]_q+\xi_{n_1+\cdots+n_s, \frac{1}{q}}.$$
 \endproclaim

\vskip 20pt

\quad Taekyun Kim

\quad Division of General Education-Mathematics, Kwangwoon
University, Seoul

\quad 139-701, S. Korea
 e-mail:\text{ tkkim$\@$kw.ac.kr}


\vskip 20pt


 \Refs \widestnumber\key{999999}
\ref \key 1
\by M. Acikgoz, S. Araci \paper  A study on the integral of the product of several type
Bernstein polynomials \jour  IST Transaction of Applied
Mathematics-Modelling and Simulation, \yr 2010\endref

\ref \key 2 \by T. Kim\paper A note on $q$-Bernstein polynomials\jour Russ. J. Math. Phys( accepted) \vol \yr  \pages \endref

\ref\key3 \by T. Kim \paper Barnes type multiple $q$-zeta function and $q$-Euler
polynomials \jour J. Phys. A: Math. Theor. \vol 43 \yr 2010 \pages 255201,
11pp\endref

\ref\key4 \by T. Kim \paper   New approach to $q$-Euler polynomials of higher order \jour Russ. J.  Math.
Phys. \vol  17 \yr 2010\pages 218-225\endref

\ref \key 5 \by T. Kim \paper  $q$-Euler numbers and polynomials associated
with $p$-adic $q$-integral \jour J. Nonlinear Math. Phys. \vol 14\yr 2007 \pages 15-27\endref

\ref\key 6
\by V. Gupta, T. Kim, J. Choi, Y.-H. Kim \paper Generating function for
$q$-Bernstein, $q$-Meyer-K\"{o}nig-Zeller and $q$-Beta basis \jour
Automation Computers Applied Mathematics \vol 19 \yr 2010 \pages 7-11\endref

\ref\key 7 \by  L.C. Jang, W.-J. Kim, Y. Simsek \paper A study on the $p$-adic integral representation on $\Bbb Z_p$ associated with Bernstein and Bernoulli polynomials \jour Advances in Difference Equations\vol 2010 \yr 2010 \pages Article ID 163217, 6pp
\endref

\ref\key 8\by   Y. Simsek, M. Acikgoz \paper A new generating function of
$q$-Bernstein-type polynomials and their interpolation function \jour Abstract and Applied Analysis
 \yr 2010 \vol 2010\pages Article ID 769095, 12pp\endref

\ref\key 9 \by T. Kim, L. C. Jang, H. Yi \paper A note on the modified $q$-Bernstein polynomials \jour Discrete Dynamics in Nature and society
 \vol 2010 \yr 2010 \pages Article ID 706483, 13pp \endref

\ref\key10 \by   B. A. Kupershmidt \paper Reflection symmetries  of $q$-Bernoulli
polynomials \jour J. Nonlinear  Math. Phys. \vol 12 \yr 2005 \pages 412-422\endref

\ref\key11 \by  T. Kim \paper Some identities on the $q$-Euler polynomials of higher order and $q$-stirling numbers by the fermionic $p$-adic
integrals on $\Bbb Z_p$ \jour Russ. J. Math. Phys.
\vol 16\yr 2009 \pages 484-491\endref










\endRefs
\enddocument